\documentclass[11pt]{article}
\usepackage[all]{xy}
\usepackage{lscape}
\usepackage[ansinew]{inputenc}
\usepackage[psamsfonts]{amssymb}
\def\bk{I\!\!k} 

\newcommand{\bh}[0]{\ensuremath{\bold{H}}}
\newcommand{\be}[0]{\ensuremath{\bold{E}}}

\title{String spectral sequences}
\author{ Jean-Fran\c cois LE BORGNE }

\begin{document}

\maketitle

\begin{abstract} The loop homology of a closed orientable manifold $M$ of
dimension $m$ is the  commutative  algebra  $\mathbb
H_\ast  (LM) = H_{\ast  +m}(LM)$ equipped with the   Chas-Sullivan loop
product.
Here $LM$ denotes the free loop on $M$. There is also a homomorphism of graded
algebras
$I$ from  $\mathbb
H_\ast  (LM)$ to the Pontryagin algebra $H_\ast (\Omega M)$. The primary
purpose of
this paper is to study the loop homology and the homomorphism $I$, when $M$ is a
fibered
manifold,  by mean of the associated Serre spectral sequence. Some related
results are
discussed, in particular the link between $I$ and the Cohen-Jones-Yan
spectral sequence.
\end{abstract}

\vspace{5mm}\noindent {\bf AMS Classification} : 55P35, 54N45,55N33, 17A65,
81T30, 17B55

\vspace{2mm}\noindent {\bf Key words} : free loop space, loop homology,
Serre spectral sequence.

\vspace{1cm}

\centerline{\sc  Introduction.}

\vspace{3mm}
In this paper  $\bk$ is   fixed commutative ring, (co)chain complexes,
(co)homology are
with coefficients in $\bk$.

Let
$M$ be a 1-connected closed  oriented
$m$-manifold and let
$LM $ be its loop space. Chas and Sullivan \cite{CS} have constructed  a  natural
product on   the desuspension of the homology of the free loop space
$
\mathbb H_\ast(LM):= H_{\ast+m}(LM)$ so that $\mathbb H_\ast(LM)$ is a
commutative
graded algebra. This product is  called the {\it loop product}.  The purpose of
this
paper is to compute this algebra for $m$-manifolds which appear as the total
space of a
fiber bundle,  (for example Stiefel manifolds). 
To be more precise let us recall that  if   $X$ and $Y$ are two Hilbert
connected  smooth oriented   manifolds without boundary and if  $ f: X \to Y $ is  a
smooth orientation preserving  embedding of codimension $k< \infty$ there is a
well
defined   homomorphism of
$H_\ast(Y)$-comodules
$$
f_! : H_\ast(Y)  {\to} H_{\ast-k} (X)\,,
$$
called the {\it homology shriek map} $ f_!$ (see section 2 for details).

Our main result, consists to  show that if $ f : X \to Y$ is  a  {\it fiber
embedding} (definition below)  then  $f_!$ behave really  nicely  with associated
Serre
spectral sequences. The remaining results of the paper are consequences of our
main
result.

A fiber embedding $(f,f^B)$ is a commutative diagram
$$
\begin{array}{cccc}
X&\stackrel{f} \to& X'\\
p\downarrow && \downarrow p'\\
B &\stackrel{f^B}\to& B'
\end{array}
$$
where
$$
(\ast)\quad
\left\{
\begin{array}{lll}
a) & X, X',B \mbox{ and } B' \mbox{ are connected Hilbert  manifolds
without
boundary}\\
b) & f \mbox{ (resp. } f^B\mbox{) is a  smooth embedding of
finite codimension } k_X
\mbox{(resp. } k_B\mbox{)}\\
c)&p \mbox{ and } p' \mbox{ are Serre fibrations}\\
d)&\mbox{for some } b \in B \mbox{ the induced map}\\
&\hspace{32 mm} f^F : F:=p^{-1}(b)\to {p'}^{-1}(f(b)):=F'\\
&\mbox{ is an embedding of finite codimension } k_F \\
e)&\mbox{embeddings } f, f^B  \mbox{ and } f^F
\mbox{admits Thom classes}
\,.
\end{array}
\right.
$$

\vspace{5mm}
\noindent{\bf First part of the main result.} {\it  Let $f: X \to X'$ be a fiber embedding as
above. For each $n\geq 0 $ there
exist filtrations
$$
\begin{array}{lcl}
&\{0\}\subset F_0C_n(X)\subset F_1C_n(X)\subset ... \subset F_nC_n(X)=C_n(X)\\

&\{0\}\subset F_0C_n(X')\subset F_1C_n(X')\subset ... \subset F_nC_n(X')=C_n(X')
\end{array}
$$
and a chain representative  $
f_! : C_\ast(X')  {\to} C_{\ast-k_X} (X)\
$
of $
f_! : H_\ast(X')  {\to} H_{\ast-k_X} (X)
$ satisfying:
$$
f_!\left( F_\ast C_\ast (X')\right) \subset   F_{\ast-k_B} C_{\ast-k_X} (X')\,.
$$}
The above filtrations induce the Serre spectral sequence of the fibration $p$
and $p'$,
denoted $\{E^r[p]\}_{r\geq 0}  $ and $\{E^r[p']\}_{r\geq 0}$.
The chain  map $f_!$ induces a
homomorphism of bidegree $(-k_B,-k_F)$ between the associated  spectral sequences
$$
\{ E^r(f_!)\} : \{E^r[p']\}_{r\geq 0} \to \{E^r[p]\}_{r\geq 0}\,.
$$
{\bf Second part of the main result.}

{\it There exists a chain representative  
$
f^B_! : C_\ast(B')  {\to} C_{\ast-{k_B}} (B)\mbox{ (resp. } f^F_! : C_\ast(F')  {\to}
C_{\ast-k_F} (F)\mbox{)}
$
of $ f^B_! : H_\ast(B')  {\to} H_{\ast-k_B} (B) $ (resp. of $ f^F_! : H_\ast(F')  {\to}
H_{\ast-k_F} (F) $ ) such that
$
\{ E^2(f_!)\} = H_\ast(f_{B!} ; {\cal H}_\ast(f^F_!)) : E^2_{s,t} [p'] =
H_s(B' ;
{\cal H}_t(F')) \to  E^2_{s-{k_B}  ,t-{k_F}} [p] = H_{s-{k_B}}(B ; {\cal
H}_{t-{k_F}}t(F))\,,
$
where  $ {\cal H}(-)$ denote the usual  system of local coefficients.}

\vspace{3mm}
Given a first quadrant homology spectral sequence $\{E_r\}_{r\geq 0}$  it is natural, in
our
context,  to defined the {\it  $(k_B,k_F)$-regraded spectral sequence}  $\{{\mathbb
E}^r\}_{r\geq 0}$ by ${\mathbb E}_{\ast, \ast} ^r = {E}_{\ast-k_B, \ast-k_F} ^r\,.$

Hereafter, we apply the main result in four different setting:

1. Classical intersection theory

2. Chas and Sullivan string theory concerning the free loop space $LM$.

3. The restricted Chas and Sullivan theory concerning the free loops
on $M$ based in a submanifold $N$ of $M$.

4. The path space: $M^I$.

{\bf 1. Classical intersection theory.} Recall also that if $M$ is $m$-dimensional oriented closed manifold then the
desuspended
homology of $M$ :
$$
{\mathbb H}_\ast(M)= H_{\ast +m}(M)
$$
is a commutative graded algebra for the intersection product \cite{Br}. Since the
intersection
product $x \otimes y \mapsto x \bullet y$ is defined by the composition:
$$
H_\ast (M) \otimes H_\ast(M)  \stackrel{\times }{\to} H_\ast(M\times M)  \stackrel{\Delta_!}{\to} H_{\ast -m}(M)\,,
$$
where $\Delta : M \to M\times M $ denotes the diagonal embedding and $\times$  the cross
product. We deduce from the main result:

\noindent{\bf Proposition 1} { \it  Let $N\to X \stackrel{p}\to M$ be a fiber
bundle such
that

a) $N$, (resp. $M$)  is a finite dimensional smooth closed  oriented manifold of
dimension
$n$ (resp. $m$)

b) $M$ is a connected space and   $\pi_1(M)$ acts trivially on $H_\ast(N)$.

The $(m,n)$-regraded Serre spectral sequence $\{{\mathbb E}^r[p]\}_{r\geq 0}$ is
a multiplicative  spectral sequence  which congerges to the algebra  ${\mathbb
H}_\ast (X)$
and such that the tensor product of graded algebras ${\mathbb H}_\ast (M)\otimes
{\mathbb
H}_\ast (N)$ is a subalgebra of ${\mathbb E}^2[p]$. In particular if $H_\ast (M)$
is
torsion free then
$${\mathbb E}^2[p]= {\mathbb H}_\ast (M)\otimes {\mathbb H}_\ast (N)\,.$$}
By Poincar\'e duality, we recover the multiplicative structure on the
Serre spectral sequence in cohomology for the cup produit. 

{\bf 2. Chas and Sullivan string theory.}
Chas and Sullivan give  a geometric construction of the loop product. An other
description
of the loop product is the following, \cite{CJY}. One consider the diagram
$$
\xymatrix{
LM \ar[d]^{ev(0)} & LM \times_M LM \ar[l]_{Comp}
\ar[d]^{ev_{\infty}} \ar[r]^{\tilde{\Delta}} & LM \times LM
\ar[d]^{ev(0) \times ev(0)} \\
M & M \ar[l]^= \ar[r]^{\Delta} & M \times M }
$$

where

a) $ev(0)(\gamma)= \gamma(0) \,, \gamma \in LM$,

b)  the  right hand square is a pullback diagram

c) $Comp$ denotes composition of free loops.

Here we assume that $LM$ is a Hilbert manifolfd when we restrict to ``smooth
loops''
and thus  $(\tilde{\Delta}, \Delta)$ is a fiber  embedding of codimension $m$,
\cite{Me}.

The  loop  product $x \otimes y \mapsto x \circ y$ is defined by the composition
$$    H_\ast(LM) ^{\otimes 2}\stackrel{\times}{\to}
H_\ast(LM^{\times 2})\stackrel{\tilde{\Delta_!}}{\to}  H_{\ast-m} (LM \times _M LM
)
\stackrel{H_\ast(Comp)}{\to} H_{\ast -m} (LM)\,.
$$
Here we assume that $LM$ is a Hilbert manifolfd when we restrict to ``smooth
loops''
and thus $\tilde{\Delta}$ is a fiber  embedding of codimension $m$,
\cite{Me}.

The desuspended homology of $LM$,
$
{\mathbb H}_\ast (LM)=H_{\ast +m}(LM)
$
is a commutative graded algebra \cite{CS}.

\noindent {\bf Theorem A: The Cohen-Jones-Yan spectral sequence.} {
  \it Let $M$ be a smooth closed oriented $m$-manifold. The
$(m,0)$-regraded Serre spectral sequence, $\{{\mathbb
E}^r[ev(0)]\}_{r\geq 0}$, of the loop fibration,  $$
\Omega M \stackrel{i_0} \to LM \stackrel{ev(0)}\to  M $$
is a multiplicative spectral sequence which congerges to the  graded algebra
${\mathbb
H}_\ast(LM)= H_{*+m}(LM)$ given at the $E^2$-level by $\mathbb{E}^2[ev(0)] = {\mathbb H}_\ast (M; {\cal H}_\ast( \Omega M))$ when ${\cal H}_\ast(
\Omega
M))$ denotes the usual system of  local coefficients.

Furthermore, if we suppose that $\pi_1(M)$ acts trivially on $\Omega M$
then the tensor product of the graded algebra
${\mathbb H}_\ast(M)$ with  the Pontrjagin algebra $H_\ast (\Omega M)$ is a
subalgebra
of the
${\mathbb
E}^2[ev(0)] = {\mathbb H}_\ast (M; H_\ast( \Omega M))$.}

\vspace{3mm}

\noindent {\bf Theorem B: The string Serre spectral sequence.} {   \it
  Under hypothesis of  Proposition 1, the
$(m,n)$-regraded
Serre spectral sequence $\{{\mathbb E}^r[Lp]\}_{r\geq 0}$ of the Serre fibration
$$
LN \stackrel{Li}\to LX \stackrel{Lp} \to LM$$
is a multiplicative spectral sequence which converges  to the algebra ${\mathbb
H}_{\ast}(LX)$ and such that  the tensor product of graded algebras ${\mathbb
H}_\ast (LM)\otimes {\mathbb H}_\ast (LN)$ is a subalgebra of ${\mathbb
E}^2[Lp]$.
In particular if $H_\ast (LM)$ is torsion free then
$${\mathbb E}^2[p]= {\mathbb H}_\ast (LM)\otimes {\mathbb H}_\ast
(LN)\,.$$
We call this spectral sequence the string ($(m,n)$-regraded) Serre
spectral sequence.}

\vspace{3mm}

Theorem A is a   generalization to the non 1-connected case  of the result obtained
in
\cite{CJY}. Spectral sequence in theorem B is naturally related  to the spectral sequence
considered in Proposition 1 by mean of each
of
the following diagrams
$$
\begin{array}{cccccc}
LN &\stackrel{Li}\to& LX &\stackrel{Lp} \to& LM\\
ev(0)^N\downarrow && ev(0)^X \downarrow  && \downarrow   ev(0)^M\\
N &\stackrel{Li}\to&X &\stackrel{Lp} \to&M\end{array}
\,, \quad
\begin{array}{cccccc}
LN &\stackrel{Li}\to& LX &\stackrel{Lp} \to& LM\\
\sigma ^N\uparrow && \sigma ^X\uparrow && \uparrow\sigma ^M\\
N &\stackrel{Li}\to&X &\stackrel{Lp} \to&M\end{array}
$$
where $\sigma ^X$ denotes the canonical section of $ev(0)^x$. Indeed both $ev(0)^X$
and
$\sigma ^X$ induces  homomorphisms of graded algebras between  ${\mathbb H}_\ast
(X)$
and   ${\mathbb H}_\ast (LX)$.

{\bf 3. Restricted Chas and Sullivan theory and intersection morphism.}
Let $i:N \hookrightarrow M$ be a smooth finite codimensionnal
embedding and $\tilde{i}:L_NM = \{\gamma \in LM ; \gamma(0)
\in N \} \hookrightarrow LM$ be the natural inclusion. Under
convenient hypothesis (see paragraph 4), $\mathbb{H}_*(L_NM)=H_{*+n}(L_NM)$
is a graded commutative algebra and Theorems A and B translate in
Theorems A' and B' as stated and proved in paragraph 4.

The case $N = pt$ is particulary interresting since $L_{pt}M = \Omega M$ and the restricted homomorphism
$\tilde{i}_!:\mathbb{H}_*(LM) \to H_*(\Omega M)$ is by definition the
{\it intersection morphism} denoted $I$ in \cite{CS}.
We prove:

{\bf Proposition 3}
{\it Let $a \in \mathbb{H}_{-m}(LM)$ be the homology class representing a
  point in $\mathbb{H}_*(LM)$ and $\mu_a: \mathbb{H}_{*}(LM) \to \mathbb{H}_{*-m}(LM)$  $x \mapsto
  a \circ x $ be multiplication by $a$.
Then, $\tilde{i}_* \circ \tilde{i}_! = \tilde{i}_* \circ I = \mu_a$.}

{\bf Theorem C}
{\it Assume that for each $k \geq 0$, $H_k(\Omega M)$ is finitly
  generated. Then, the three following propositions are equivalent:

(a)  $I$ is onto

(b)  The differentials $\{d_n\}_{n \geq 0}$ of the  Cohen-Jones-Yan spectral sequence
vanish for $n \geq 2$ so that the spectral sequence collapses at the $E^2$-level

(c) $\tilde{i}_*$ is injective}

Explicit computations are given in subsection 4.6.

{\bf 4. The path space $M^I$.}
Consider the fibration
$
\Omega M \stackrel{i_1}\to M^I \stackrel{ev(0), ev(1)} \to  M\times M \,.
$
First, one should observe   that the construction of the loop product extends to
$M^I$ so that ${\mathbb H}_\ast(M^{I})$ is a graded commutative algebra
isomorphic to
${\mathbb H}_\ast(M)$. Secondly, there exists on ${\bf H}_{\ast} (M\times M)=
{\mathbb H}_{\ast-m} (M\times M)$ a natural structure of associative
algebra (not
commutative
and without unit). With these two datas, we prove Theorem A''
which extends Theorem A (See the section 5 for more details).

{\bf Acknowledgement}

We would like to thank David Chataur, Jean-Claude Thomas and Luc
Mennichi for
helpful comments that have considerably improved the paper.

\vspace{3mm}

The paper will be organized as follows. In Sect. 1 we will recall the  definition
of the shriek map associated to a smooth embedding. In Sect. 2 we will prove the
main result. In Sect. 3 we will prove Proposition 1, Theorem A and B. In
Sect. 4 we will state and prove Theorems A' and B' Proposition 2 and 3, prove Theorem C and give some examples. In Sect. 5 we will state
and  prove Theorems A'' and Proposition 4 and give an example.

\vspace{3mm}
\centerline{\sc 1. Recollection on Thom-Pontrjagin theory.}

The purpose of this section is to precise the definition of the shriek
map of an embedding \cite{CT}.
 
{\bf 1.1} Let us recall that the Thom class $\tau _{E,E'}$ of   a disk
bundle
$$
(p) \qquad   (D^{k},S^{k-1}) \to (E,E') \stackrel{p}\to B
$$
is an element $\tau_{p }\in H^k(E,E')$ whose restriction to each fiber $(D^k,
S^{k-1})$ is a generator of $H^k(D^k, S^{k-1})$. Every disk bundle has a Thom
class if
$\bk = {\mathbb F}_2$. Every oriented disk bundle has a Thom class with any ring
of
coefficients $\bk$. If the  disk bundle $(p)$  has a Thom class $\tau$  then the
composition
$$
H_{\ast }(E,E')\stackrel{\tau _p\cap -}{\to }H_{\ast -d}(E) \stackrel{
H_\ast(p)}{\to}    H_{\ast - d}(B)
$$
is an isomorphism of graded modules, called the {\it homology Thom isomorphism of
the fiber  bundle $(p)$}.

{\bf 1.2} 
Let $N$ and $M$ be two smooth closed  oriented manifolds and  $ f: N
\to M
$ be a smooth orientation preserving  embedding. For simplicity, we  identify
$f(N)$ with
$N$. Assume further, that $N$ is closed subspace of $M$ and that $M$ admits a
partition
of unity  then  there is a splitting
$$
T(M)_{|_N} = T(N) \oplus \nu_f
$$
where $\nu_f$ is called the normal fiber bundle of the embedding. The isomorphism
class
of $\nu_f $ is well defined and depends only of the isotopy class of $f$.

Assume that the rank of $\nu_f=k$ and denote
by
$$(D^k ,S^{k-1}) \to (D^k_{p_f}, S^{k-1}_{p_f}) \stackrel{p_f}{\to} N $$
the associated disk bundle. The restriction of the exponential
map induces an isomorphism $\Theta$ from $(D^k_{p_f}, S^{k-1}_{p_f})$ onto the tubular neighborhood
$(\mbox{Tube}\, f\,,
\partial \mbox{Tube}\, f)$.

The Thom class of the embedding $f$ is the Thom class, whenether it exists,  of
the disk
bundle $p_f$.

The definition of $f_!:H_*(M) \longrightarrow H_{*-k}(N)$ is given by
the following composition:
$$
\xymatrix{
H_*(M) \ar[r]^{j^M} & H_\ast(M,M-N) \ar[r]^{Exc} & H_\ast
(\mbox{Tube}\, f ,\partial \mbox{Tube}\, f) \ar[r]^{\Theta_*} &
H_\ast((D^k_{p_f}, S^{k-1}_{p_f})) \ar[r]^{p_*} & H_*(N) }
$$
where $Exc$ denotes the excision isomorphism,  $j^B :   H_\ast(B) \to H_\ast
(B,A) $
the canonical  homorphism   induced by the inclusion $A
\subset B$.
\vspace{3mm}

\centerline{\sc 2. Proof of the main result.}

The aim  of this section is  to prove our main result. For this
purpose following \cite{CJY} we describe the homology shriek map at the chain level and the Serre filtration.

{\bf 2.1 Shriek map of an embedding at the chain level.}
We use the same notations as in section 1 and we define
$f_!:C_*(N) \to C_{*-k}(M)$.
Let $t:M \to Tube f / \partial Tube f$ the Thom collapse map wich is
identity on the interior of the tubular neighbourood $Tube f$ and
collapses the rest on a single point. By construction, this map is
continuous. We denote by $t_{\sharp}$ the map induced at the chain
level by $t$.
Chose $\tilde{\tau}_{p_f} \in  C^k(D^k_{p_f}, S^{k-1}_{p_f})$
representing $\tau_{p_f}$ and let $f_!$ at the chain level be the composition:
$$ 
C_*(M) \stackrel{t_{\sharp}}{\to}  C_*(Tube f / \partial Tube f)   \stackrel{can}{\to}  C_*(Tube f , \partial Tube f)   \stackrel{\Theta_{\sharp}}{\to}
  C_*(D^k_{p_f}, S^{k-1}_{p_f}) 
  \stackrel{\psi}{\to} C_{*-k}(N)
$$
where $can$ is the algebraic natural application $can:C_*(A/B) \to
C_*(A, B)$, $\Theta_{\sharp}$ the chain map induced by the exponential
map and $\psi$ defined as follow:
since the cap product is well defined at the level of the cochains,
we define:
$$\psi:C_*(D^k_{p_f}, S^{k-1}_{p_f}) \longrightarrow   C_{*-k}(N) \qquad
\alpha \mapsto (p_f)_{\sharp}(\tilde{\tau}_{p_f} \cap \alpha).$$
The definition of $\psi$ depends of the choice of
$\tilde{\tau}_{p_f}$ but in homology, each choice induces the same morphism: the
Thom isomorphism.

{\bf 2.2 Serre filtration.}
Following \cite{MC} there is a filtration of $C_*(X)$.
$F_{-1}=0 \subset F_0 \subset F_1 \subset ... \subset F_n \subset
... \quad $ defined as follows:
for $p \geq q$, define $(i_0,i_1,...,i_p):\Delta^q \longrightarrow \Delta^p$ the
linear map wich maps the vertex $v_k$ of the standard $q$-simplex
$\Delta^q$ to the vertex $v_{i_k}$ of  $\Delta^p$ and put 
$
F_pS_q(X)=\{\sigma \in S_*X / 
\mbox{ there exists } 
\Sigma \in S_*(B)  \mbox{ such that }
p \circ \sigma = \Sigma \circ (i_0,i_1,...,i_q) \}$.
The linear extension of this filtration provides a filtration of
$C_*(X)$.
This filtration leads to the construction of the Serre
spectral sequence.
To continue the proof of the main result, we need two lemmas. 

{\bf 2.3 Lemma 1}
{\it A morphism of fibrations

$$
\xymatrix{
F \ar[d] \ar[r]^{i\vert_F}& F' \ar[d] \\
 E   \ar[d]^p \ar[r]^{i} & E' \ar[d]^{p'}\\ 
B \ar[r]^{\tilde{i}}  & B' } 
$$
such that $\tilde{i}(B)$ and $p'(E')$
intersect transversally in $B'$ factorizes as 
$$
\xymatrix{
 F   \ar[r]^{i\vert_F} \ar[d]& F'  \ar[r]^{=} \ar[d]&  F' \ar[d]  \\
 E \ar[r]^{i_{2E}} \ar[d]^p &   E \times_{B'} B  \ar[r]^{{i_{1E}}} \ar[d]^{p''} & E' \ar[d]^{p'}\\
 B  \ar[r]^{=}  & B \ar[r]^{\tilde{i}} & B'}
$$
where $p''$ is the pull-back fibration of $p'$ along $i$.
Moreover, the upper left square is a pull-back diagram.}

Proof:
Consider the commutative diagram:
$$
\xymatrix{
 F   \ar@{-->}[r]^{i_{\vert_F}} \ar[d]^{j_1}& F'  \ar[r]^{=} \ar[d]^{j_2} &  F' \ar[d]  \\
 E \ar@{.>}[r] \ar[rd]^p &   E \times_{B'} B  \ar[r] \ar[d] & E' \ar[d]^{p'}\\
  & B \ar[r]^{\tilde{i}} & B'}
$$
where the right hand part is a pull-back diagram and the dotted arrow is obtained by universal property.
We want to show that the left hand square is a pullback diagram
where $j_1$ and $j_2$ are inclusions of fibers in the total space.
For this purpose we consider the following diagram
$$
\xymatrix{
F \ar[r]^{i{\vert}_F} \ar[d]^{j_1} & F' \ar[d]^{j_2} \ar[r] & pt \ar[d]\\
E \ar[r] & E' \times_{B'} B \ar[r]^{p''} & B}
$$
where $pt$ is a point of $B$. The big square is a pull-back
diagram. Therefore $F_1 \cong E \times_B pt \cong E \times_{ (E'
  \times_{B'} B)} ( (E' \times_{B'} B) \times_B pt)$
($\cong $ means homeomorphic). The right square is a pull-back thus 
$F' \cong  (E' \times_{B'} B) \times_B pt$
and $F \cong E \times_{ E' \times_{B'} B} F'.$

\rightline{$\square$}

{\bf 2.4 Lemma 2}
{\it Denote by $\nu_f$ the normal bundle of the embedding $f:X
  \hookrightarrow X'$ of the fiber embedding. If we chose $\tilde{\tau} \in
C^*(D^{k_X}(\nu_f),S^{k_{X}-1}(\nu_f))$ such that $\tilde{\tau}$ vanishes on
degenerate simplices, then  
$$
F_pC_{p+q}(D^{k_X}(\nu_f),S^{k_X-1}(\nu_f))
\stackrel{\tilde{\tau} \cap -} \longrightarrow F_{p-k_B}C_{p+q-k_X}(D^{k_X}(\nu_f),S^{k_X-1}(\nu_f))
$$
for $p$,$q$, integers.}

Proof:
Applying Lemma 1 to our fiber embedding $(f,f_B)$, we obtain the next
diagram:
$$
\xymatrix{
 F   \ar[r]^{f^F} \ar[d]& F'  \ar[r]^{=} \ar[d]&  F' \ar[d]  \\
 X \ar[r]^{f_2} \ar[d]^p &   X \times_{B'} B  \ar[r]^{f_1} \ar[d]^{p''} & X' \ar[d]^{p'}\\
 B  \ar[r]^{=}  & B \ar[r]^{f^B} & B'}.
$$

We begin by proving that right and left hand part of this diagram are fibers embedding.
Since $f^B(B)$ and $p'(X')$
intersect transversally in $B'$ \cite{K}, $X \times_{B'} B$ is a
manifold so that the pull-back
diagram of the right is a fiber embedding. The same holds for the
left part of the diagram.

Now, we will prove lemma 2 for each fiber embedding $(f_1,f^B)$ and $(f_2,id_B)$. 

1) Case  $(f_1,f^B)$.
Let $\tau_{f^B}$ and $\tau_{f_1}$ be 
respectively the Thom class of $f^B$ and $f_1$. By naturality, we have 
$\tau_{f_1}= p'^*(\tau_{f^B})$. At the chain level we
can choose a cocycle representing the Thom classes, also called
$\tau_{f^B}$ and $\tau_{f_1}$, and we have $\tau_{f_1} = p'^{\sharp}(\tau_{f^B})$.

For some $\sigma \in F_pC_{p+q}(X')$, by definition of the Serre
filtration, there exists $ \Sigma \in C_p( B')$ and $ (i_0,...,i_{p+q})$ such that $p'(\sigma)=\Sigma(i_0,...,i_{p+q})$.
Thus $p'_{\sharp} (\tau_{f_1} \cap \sigma) = p'_{\sharp} (
p'^{\sharp}(\tau_{f_B}) \cap \sigma)
=\tau_{f_B} \cap \Sigma(i_0,...,i_{p+q})
=\tau_{f_B}(i_{p+q-k_B},...,i_{p+q})  \Sigma(i_0,...,i_{p+q-k_B})$.

Since  $\tau_{f_B}(i_{p+q-k_B},...,i_{p+q}) \neq 0 $ ,   $ i_{p+q-k_B},...,i_{p+q} $
are pairwise distinct then we can write $\Sigma(i_0,...,i_{p+q-b}) $
as an element of $ C_{p-k_B}( B')$
so that $f_{1!}(\sigma) \in  F_{p-k_B}C_{p+q-k_B}( X' \times_{B'}
B) $.
  
It remains to show that $f_{1!}$ preserves the differentials (the
context being sufficiently clear to know about wich differential we
refer, all differentials are denoted by $d$).
Let  $c \in  C_*(X')$, since $d(\tau_{f_{1}})=0$ we have:
$d(\tau_{f_1} \cap c) = d(\tau_{f_1}) \cap c +
(-1)^{k_B} \tau_{f_1} \cap dc 
$$= (-1)^{k_B} \tau_{f_1} \cap dc$
This implies that  
$d(f_{1!}(c))=(-1)^{k_B}i_!(dc)$.

2) Case $(f_2,id_B)$.
Let $\omega \in  F_{p-k_B}(C_{p+q-k_B}(X' \times_{B'} B)) $ .
Then there exists $ \Omega \in C_{p-k_B}(B)$ such that $p''_{\sharp}(\omega)=\Omega(i_0,...,i_{p+q-k_B})$.
Let $\tau_{f_2} \in C^{k_F}(X \times_B' B)$ be a cochain representing the Thom class of the embedding $f_2$.
 Thus  $\tau_{f_2} \cap \omega $ is a subchain of degree 
$p+q-k_X$ thus  $\tau_{f_2} \cap \omega $ factorizes by $\Omega$. 
Then $f_{2!}(\omega)$ lies in  $F_{p-k_B}C_{p+q-k_X}(X)$.
It is now easy to complete the proof.
\rightline{$\square$}
\vspace{3mm}
{\bf 2.5 End of the proof of the main result.}
Since the Serre filtration of $C_*(X)$ is natural with respect to
fiberwise maps the three first applications defining $f_!$ as in {\bf
  2.1} induce
morphisms of differential gradued filtred module (dgfm for short).
Lemma 2 proves that the last application is a morphism of dgfm then
$f_!$ induces a morphism of dgfm of bidegree $(k_B, k_X)$.

The second part of the main result follows by classical theory of
spectral sequences \cite{MC}.

\rightline{$\square$}

\vspace{3mm}
\centerline{\sc 3. Proof of Proposition 1, Theorems A,B.}
{\bf 3.1 Proof of Proposition 1} 
For a topological space $Y$, denote by $\Delta_Y$ the diagonal embedding
of $Y$ in $Y \times Y$.
The diagonal map $\Delta_X : X \longrightarrow X \times X$ factorizes
so that we obtain the following commutative diagram :
$$
\xymatrix{
 N \times N   \ar[r] &   X \times X  \ar[r] &  M \times M 
 \\
 N \times N   \ar[r] \ar[u]^{id} &   X \times_M X   \ar[r] \ar[u]^{i} &  M \ar[u]^{\Delta_M}\\
 N   \ar[r] \ar[u]^{\Delta_N} &   X   \ar[r] \ar[u]^{\Delta'_X} &  M \ar[u]^{id}}
$$
with $i \circ \Delta'_X = \Delta_X$.
In particular $(i,\Delta_M)$ and $(\Delta'_X , id)$ are
fiber embeddings.
We apply the main result to the shriek map $\Delta_{X!}$.
Since the composition 
$$\xymatrix{
H_*(X) \otimes H_*(X) \ar[r]^{\times} & H_*(X \times X)
\ar[r]^{\Delta_{X!}} &  H_{*-{(m+n)}}(X) }$$
is the intersection product one deduces that the regraded spectral
sequence is a multiplicative spectral sequence with respect to the
intersection product. 
The $E^2$-term is given by $\mathbb{E}^2_{p,q}=\mathbb{H}_p(M ;
\mathbb{H}_q(N)) $ (since  $\pi_1(M)$ is
assumed to act trivially on $H_*(N)$ the coefficients are constant).
The naturality of the cross product provides a morphism of spectral
sequence given at the  $E^2$-term by 
$$\xymatrix{
\mathbb{H}_p(M; \mathbb{H}_q(N)) \otimes \mathbb{H}_{p'}(M;
\mathbb{H}_{q'}(N)) \ar[r]^{\times} & \mathbb{H}_{p+p'+m}(M \times M ;
\mathbb{H}_{q+q'+n}(N \times N))}
$$   
Then, $\Delta_{X!}$ induces a morphism of spectral sequence given at
the  $E^2$-term by 
$$\xymatrix{
\mathbb{H}_{p+p'+m}(M \times M ; \mathbb{H}_{q+q'+n}(N \times N))
\ar[rr]^{E^2(\Delta_{X!})} & & \mathbb{H}_{p+p'}(M;
\mathbb{H}_{q+q'}(N))}
$$
such that $E^2(\Delta_{X!}) = H(\Delta_{M!}; \Delta_{N!})$.
As a consequence, we find that $\mathbb{H}_*(M) \otimes
\mathbb{H}_*(N) \hookrightarrow \mathbb{E}^2_{*,*}$ as subalgebra
where  $\mathbb{H}_*(M) \otimes \mathbb{H}_*(N)$ is the tensor product
of algebra for the intersection product.

\rightline{$\square$}

{\bf 3.2 Proof of theorem A: the Cohen-Jones-Yan spectral sequence. }
In this section the results of \cite{CJY} are revisited and sligtly
extended. We construct the following commutative diagram:
$$
\xymatrix{
\Omega M \times \Omega M \ar[r]^{i \times i}  & LM \times LM
\ar[r]^{ev(0) \times ev(0)} & M \times M \\
\Omega M \times \Omega M \ar[r]^{i \times i} \ar[u]^{id} \ar[d]^{Comp} &  LM
\times_M LM \ar[u]^{\tilde{\Delta}} \ar[r]^{ev_{\infty}} \ar[d]^{Comp}& M
\ar[u]^{\Delta} \ar[d]^{id}\\
\Omega M \ar[r]^{i}  & LM \ar[r]^{ev(0)} & M}
$$
where the map defined as in the introduction.
Thus $(\tilde{\Delta}, \Delta)$ is a fiber embedding of codimension $m$.
Now, composition of the maps
$$
C_{*+d}(LM) \otimes C_{*+d}(LM) \stackrel{\times}{\to} C_{*+2d}(LM \times LM)
\stackrel{\tilde{\Delta}_!}{\to}   C_{*+d}(LM \times_M LM)
\stackrel{Comp_{\sharp}}{\to} C_{*+d}(LM)
$$
induces at the homology level Chas and Sullivan product denoted by $\mu$.
The Serre spectral sequence associated to the fibration 
$$
\xymatrix{
\Omega M \ar[r]^i & LM \ar[r]^{ev(0)} & M}
$$
satisfy $\mathbb{E}^2_{*,*} = \mathbb{H}_*(M; H_*(\Omega M))$ (here $\pi_1(M)$ acts trivialy on $\Omega
M$ since $M$ is arcwise connected ).

By using the main result for $(\tilde{\Delta},\Delta)$ and the
naturality of  the Serre spectral sequence for the Eilenberg-Zilbert
morphism and for $\gamma_{\sharp}$, we show extending the result of \cite{CJY} that
there is a multiplicative structure on this Serre spectral sequence
containing at the $E^2$-level the tensor product of $\mathbb{H}_*(M)$
with intersection product and $ H_*(\Omega M)$ with Pontryagin product.  
This spectral sequence of algebra converges to $\mathbb{H}_*(LM)$.

{\bf 3.3 Proof of theorem B: String Serre spectral sequence.}
For any topological space $T$, the map
$map(T,E) \to map(T,B)$ is a fibration with fiber $map(T,F)$.
Since $LT = map(S^1, T)$ and $LT \times_T LT = map(S^1 \vee S^1, T)$, 
we have the following fibrations: $\xymatrix{
LN \ar[r]^{Li} & LX \ar[r]^{Lp} & LM}$
and $\xymatrix{
LN \times_N LF\ar[r]^{Li} & LX \times_X LX \ar[r]^{Lp} & LM \times_M LM }$
Consider the pull-back diagrams
$$
\xymatrix{
X \times_M X \ar[r] \ar[d] & X \times X \ar[d]^{p \times p} \\
M \ar[r]^{\Delta_M} & M \times M }
\mbox{ and }
\xymatrix{
LX \times_{X\times_MX} LX \ar[r]^{\Delta_1} \ar[d]^{Lp
  \times_{X\times_MX} Lp} & LX \times LX
\ar[d]^{Lp \times Lp} \\
LM \times_M LM  \ar[r]^{\Delta_M} & LM \times LM. } 
$$
with
$ LX \times_{X\times_MX} LX = \{(\gamma_1,\gamma_2) \in LX \times LX \quad
; \quad p(\gamma_1(0))=p(\gamma_2(0)) \}$.
We deduce the following commutative diagram:
$$
\xymatrix{
 LN \times LN   \ar[r] &   LX \times LX  \ar[r] &  LM \times LM 
 \\
 LN \times LN   \ar[r] \ar[u]^{id} &   LX \times_{X\times_MX} LX  \ar[r] \ar[u]^{i} &  LM \times_M LM  \ar[u]^{L\Delta_M}\\
 LN \times_N LN   \ar[r] \ar[d]^{Comp_N} \ar[u]^{L\Delta_N}&   LX
   \times_X LX  \ar[r] \ar[d]^{Comp_X} \ar[u]^{L\Delta'_X} &  LM \times_M LM \ar[u]^{id} \ar[d]^{Comp_M} \\
 LN \ar[r] &   LX \ar[r] &  LM }
$$
with obviously defined map. From the definition of $\mu_X$, naturality
and the main result, we construct a product on the Serre spectral
sequence associated to the fibration 
$
\xymatrix{
LN \ar[r]^{Lj} & LX \ar[r]^{Lp} & LM}
$
To achieve the proof, we need:

{\bf Lemma} {\it 
Let $\xymatrix{F \ar[r]^j & E \ar[r]^p & B}$ be a fibration such that $\pi_1(B)$
acts trivially on $H_*(F)$, then $\pi_1(LB)$ acts trivially on $H_*(LF)$.}

Proof of the lemma: For a fixed loop $\gamma \in \Omega B$, the holonomy operation defines maps  $$\Psi_{\gamma}:F \longrightarrow F
\qquad x \mapsto \gamma.x.$$
Since $\pi_1(B)$ acts trivially on $H_*(F)$, there exists a homotopy $H_{\gamma}:F \times I \longrightarrow F \qquad
(x,t) \mapsto H_{\gamma }(x,t)$ such that $H_{\gamma}(x,0)=x$, and
$H_{\gamma}(x,1)=\gamma.x$.
For a fixed $\Gamma \in \Omega LB$, the holonomy
action of $\Omega LB$ on $LF$ (associated to the fibration $
\xymatrix{
LF \ar[r]^{Lj} & LE \ar[r]^{Lp} & LB}
$)
yields
$$\Phi_{\gamma}: LF \longrightarrow LF \qquad f \mapsto (s \mapsto
\Psi_{\Gamma(-,s)}(f(s)))$$
Now the homotopy
$$\mathcal{H}_{\Gamma}: LF \times I \longrightarrow LF \qquad f(-),t
\mapsto \mathcal{H}_{\Gamma}(f(-),t)=(s \mapsto
\mathcal{H}_{\Gamma}(-,s)(f(s),t))$$
satisfies
$$\mathcal{H}_{\Gamma}(f,0) = (s \mapsto H_{\Gamma(-,s)}(f(s),0))=f(-))$$
$$\mathcal{H}_{\Gamma}(f,1) = (s \mapsto
H_{\Gamma(-,s)}(f(s),1))=\Gamma(-,-).f(-) = \Phi_{\Gamma.f}$$ $s \in S^1$.

\rightline{$\square$}

The above lemma proves that the local coefficients in the
spectral sequence are constant.
Furthermore, $\mathbb{E}^2_{*,*} = H_{*+m}(LM; H_{*+n}(LN))$ contains
$\mathbb{H}_*(LM) \otimes \mathbb{H}(LN)$ as subalgebra.

\rightline{$\square$}

\centerline{\sc 4. Restricted Chas and Sullivan algebra and
  intersection morphism.}
\centerline{\sc Proof of theorem A',B',C and proposition 2, 3.}
{\bf 4.1 The restricted Chas and Sullivan loop product.}
Let $i:N \hookrightarrow M$ be a smooth finite codimensionnal
embedding of Hilbert closed manifold such that:

*) $M$ (resp. $N$) is of finite dimension $m$ (resp. $n$).

**) the embedding $i$ admit a Thom class.

Define $L_NM$ as the right corner of the pull-back diagram:
$$
\xymatrix{
L_NM \ar[d]^{ev(0)} \ar[r]^{\tilde{i}} & LM  \ar[d]^{ev(0)} \\
N \ar[r]^i & M }
$$

i.e $L_NM$ is the space of free loop spaces of $LM$ based on $N$.
Denote $\mathbb{H}_*(L_NM) := H_{*+n}(L_NM)$.
We define the restricted loop product $$\mu_N:\mathbb{H}_*(L_NM)
\otimes \mathbb{H}_*(L_NM) \to \mathbb{H}_*(L_NM)$$
$$\alpha \otimes \beta \mapsto \mu_N(\alpha \otimes \beta) := \alpha \circ_N
\beta$$
by putting $\mu_N= Comp_* \circ \Delta_{N!} \circ \times$ where we
consider the commutative diagram:
$$
\xymatrix{
L_NM \ar[d]^{ev(0)} & L_NM \times_N L_NM \ar[l]_{Comp}
\ar[d]^{ev_{\infty}} \ar[r]^{\tilde{\Delta_N}} & L_NM \times_N L_NM
\ar[d]^{ev(0) \times ev(0)} \\
N & N \ar[l]^= \ar[r]^{\Delta_N} & N \times N }
$$

{\bf Theorem A'}
{\it The (n,0)-regraded spectral sequence associated to
  the fibration $\Omega M \longrightarrow L_NM \longrightarrow N$ is
  multiplicative. Moreover, if $\pi_1(N)$ acts trivially on $\Omega
  M$, the $E^2-term$ of the spectral sequence contains
  $i_!(\mathbb{H}_*(M)) \otimes H_*(\Omega M)$ as subalgebra.}

Proof:
Starting with the following commutative diagram:
$$
\xymatrix{
\Omega M \times \Omega M \ar[r]  & L_NM \times L_NM
\ar[r]^{ev(0) \times ev(0)} & N \times N \\
\Omega M \times \Omega M \ar[r] \ar[u]^{id} \ar[d]^{Comp} &  L_NM
\times_N L_NM \ar[u]^{\tilde{\Delta_N}} \ar[r]^{ev_{\infty}} \ar[d]^{Comp}& N
\ar[u]^{\Delta_N} \ar[d]^{id}\\
\Omega M \ar[r]  & L_NM \ar[r]^{ev(0)} & N}
$$
the proof of Theorem A works as well to prove Theorem A'.

\rightline{$\square$}

{\bf Proposition 2}
{\it  $\tilde{i}_!$ is a morphism of algebra and induces
 a morphism of multiplicative spectral sequence:
$E^*( \tilde{i}_!):\mathbb{E}^*(LM) \longrightarrow
 \mathbb{E}^*(L_NM)$.}

Proof:
Observe first that $(\tilde{i},i)$ is a fiber embedding so that the
result comes immediately from the main result and from the following commutative diagram:
$$
\xymatrix{
L_NM \times L_NM \ar[r]^{\tilde{i} \times \tilde{i}} & LM \times LM \\
L_NM \times_N L_NM \ar[u] \ar[r]^{\tilde{i} \times_N \tilde{i}}
\ar[d]^{Comp} & LM \times_M LM \ar[u] \ar[d]^{Comp} \\
L_NM \ar[r]^{\tilde{i}} & LM }.
$$
Indeed $\tilde{i}_! \circ \mu = \mu_N \circ \tilde{i}_! \otimes
\tilde{i}_!$.

\rightline{$\square$}

Now, let state theorem B'.
Let $\xymatrix{N \ar[r]^{i} & X \ar[r]^p & M}$ be the fibration of
  theorem B and $j_V:V \hookrightarrow M$ an embedding satisfying conditions (*)
  and (**). Construct the induced bundle $p_{\vert Y}$ from the pull-back diagram:
$\xymatrix{
U \ar[d]_{i \vert_U} \ar[r]^{j_{\vert_U}} & N \ar[d]^i \\
Y \ar[d]^{p_{\vert_Y}} \ar[r]^j & X \ar[d]^p \\
V \ar[r]^{j_V} & M }$

{\bf Theorem B'}
{\it In the above  situation there is a morphism of multiplicative spectral
  sequence 
$\xymatrix{
\mathbb{E}^*_{*,*}(Lp) \ar[r]^{E(\tilde{j_!})} &
\mathbb{E}^*_{*,*}(p_{\vert Y})}$ 
given at the $E^2$-level by $E(\tilde{j_!}) = H_*(j_{V!};j_{\vert_U!})$.}

Proof:
The pull-back diagram:
$\xymatrix{
L_UN \ar[d]_{Li \vert_U} \ar[r]^{\tilde{j_{\vert_U}}} & LN \ar[d]^{Li} \\
L_YN \ar[d]^{Lp_{\vert_Y}} \ar[r]^{\tilde{j}} & LX \ar[d]^{Lp} \\
L_VN \ar[r]^{\tilde{j_V}} & LM }
$
is a fiber embedding. The theorem comes directly from the main result.

\rightline{$\square$}

{\bf Remark:}
If $M = N \times U$, and if $p:M \longrightarrow N $ is the projection
on the first factor and if $i: N \hookrightarrow M $ is a standard
embedding, then $p_* \circ \tilde{i}_!: \mathbb{H}_*(LM)
\longrightarrow \mathbb{H}_*(L_NM) \longrightarrow \mathbb{H}_*(LN)$
is the projection $\mathbb{H}_*(LN) \otimes \mathbb{H}_*(U)
\longrightarrow \mathbb{H}_*(N)$.

{\bf 4.3 Proof of proposition 3.}
We consider the following commutative diagram:
$$
\xymatrix{
pt \times  LM  \ar[r]^{j} & LM \times LM \\
pt \times \Omega M  \ar[r] \ar[u]^{id \times \tilde{i}} \ar[d]^{Comp}
& LM \times_M LM  \ar[u]^{\tilde{\Delta}} \ar[d]^{Comp} \\ 
\Omega M  \ar[r]^{\tilde{i}} & LM 
}
$$
We observe that $\tilde{\Delta}_{\vert_{pt \times \Omega M}} = id \times
\tilde{i}$ and that $Comp:pt \times \Omega M \to \Omega M$ is homotopic to
the identity. Denote by $EZ$ the cross product. The map $\mathbb{H}_{*}(LM) \to
\mathbb{H}_{*-m}(LM)$  $x \mapsto Comp_* \circ \tilde{\Delta}_!
\circ EZ(pt,x) $ is in fact multiplication by $a$, with $a$ the homology class of $pt$ in $\mathbb{H}_{*-m}(LM)$.
The other map $\tilde{i}_* \circ Comp_* \circ (id \times
\tilde{i}_!) EZ(pt,-)$ is equal to  $\tilde{i}_* \circ \tilde{i}_!$.

\rightline{$\square$}

{\bf 4.4 Beginning of the proof of theorem C.}
(1) Assume that for $n \geq 2$, all the differentials of the  Cohen-Jones-Yan spectral
sequence vanish, then the homomorphism $E^*_{*,*}(\tilde{i}_!)$ induced by
the fiber embedding 
$$
\xymatrix{
\Omega M    \ar@{^{(}->}[r]^{id} \ar[d]^{id} &  \Omega M \ar[d] \\
\Omega M   \ar@{^{(}->}[r]^{\tilde{i}} \ar[d]^{ev(0)}
& LM  \ar[d]^{ev(0)} \\ 
pt  \ar@{^{(}->}[r]^{i} & M }
$$ 
is clearly onto. Denote by $E^*_{*,*}(1)$ the spectral sequence
associated to the left fibration and by  $\mathbb{E}^*_{*,*}(2)$ the
$(d,0)$-regraded spectral sequence associated to the right fibration.    
Then, at the aboutment, $E^{\infty}(\tilde{i}_!)$ is onto on the graded
space of $H_*(\Omega M)$ then $I$ is onto.
\vspace{3 mm}

(2) We begin by proving that all the differentials
starting from $\mathbb{E}^*_{0,*}(2)$ vanish. We write the naturality
of the Serre spectral sequence to the shriek map of the fiber
embedding $(\tilde{i},i)$ shown in the main result. 
We have the following commutative diagram:
$$\xymatrix{
H_*(\Omega M) \simeq \mathbb{E}^2_{0,*}(2)
\ar[d]^{E^2(\tilde{i}_!)=i_! \otimes id}& \supseteq
\mathbb{E}^3_{0,*}(2) & \supseteq ...& \supseteq
\mathbb{E}^{\infty}_{0,*}(2) \ar[d]^{E^{\infty}(\tilde{i}_!)} &
\mathbb{H}_*(LM) \ar@{>>}[l] \ar[dl]^{\tilde{i}_!} \\
H_*(\Omega M) \simeq E^2_{0,*}(1) & = E^3_{0,*}(1) & =
...& = E^{\infty}_{0,*}(1) = H_*(\Omega M)  & }
$$
If $I = \tilde{i}_!$ is onto and since each $H_k(\Omega M)$ is finitely
generated, we have
$$\mathbb{E}^2_{0,*}(2) = \mathbb{E}^3_{0,*}(2) = ... =
\mathbb{E}^{\infty}_{0,*}(2). $$  Thus the differentials
starting from $\mathbb{E}^*_{0,*}(2) $ vanish.  
We remark that in this case, the morphism at the top right of the
above diagram is in fact $I$. The existence of the canonical section
$M \to LM$ implies that the differentials starting from
$\mathbb{E}^*_{*,0}$ vanish. The multiplicative structure of
$\mathbb{E}^2_{*,*}(2)$ implies that all the differentials of the
Cohen-Jones-Yan spectral sequence vanish.

To end the proof of Theorem C we need:

{\bf Lemma}
{\it If one of the differential of the Cohen-Jones-Yan spectral
  sequence is non zero, then there exists a non-zero differential
  arriving on $\mathbb{E}^*_{-d,*}(2)$.}

Proof :
Denote by $[M] \in \mathbb{H}_*(M)$ the fundamental class of $M$ and
by $1_{\omega}$ the unit of $H_*(\Omega M)$. 
Assume that there exist a non-zero differential in the
Cohen-Jones-Yan spectral sequence. We consider the first page of the
spectral sequence where there is a non-zero differential. This
page is isomorphic to $\mathbb{E}^2_{*,*}$ as an algebra. Since the Cohen-Jones-Yan spectral
  sequence is multiplicative, there exists a non-zero differential
  starting from a generator of $\mathbb{E}^2_{*,*}$ namely an element
$[M] \otimes \omega$  of $\mathbb{E}^2_{0,*}$.
Let $x \otimes \omega' = d([M] \otimes \omega)$ and $y \in
\mathbb{H}_*(M)$ such that $x \bullet y = *$ with $* \in
\mathbb{H}_{-d}(M)$ representing a fixed point.
Then, $d(y \otimes \omega)=\pm d(y \otimes 1_{\omega} \circ [M] \otimes
\omega)=\pm  y \otimes 1_{\omega} \circ d([M] \otimes
\omega)= \pm  y \otimes 1_{\omega} \circ x \otimes \omega' = \pm x
\bullet y \otimes \omega' = \pm * \otimes \omega'$. Moreover $d(y \otimes 1_{\omega})=0$
because of the existence of a section.

\rightline{$\square$}

{\bf 4.5 End of proof of theorem C.}

{\bf Observation}
If $f:M \to N$ is a map between Poincar\'e duality manifolds,
then $f_!$ is onto iff $f_*$ injective.
We prove that this result is true for the embedding $\tilde{i}: \Omega M \hookrightarrow LM$.

Proof:
(1) Assume that for $n \geq 2$, all the differentials of the  Cohen-Jones-Yan spectral
sequence vanish, then the homomorphism $E^*_{*,*}(\tilde{i}_*)$ is
clearly injective. Furthermore, we have the injective map $\mathbb{E}^{\infty}_{-d,*}(2) \hookrightarrow \mathbb{H}_*(LM)$.
The composition of this two application is $\tilde{i}_*$. 

\vspace{3mm}
(2) Assume $\tilde{i}_*$ is injective.
From the naturality of the Serre spectral sequence for $\tilde{i}_*$,
we deduce the following commutative diagram:
$$
\xymatrix{
H_*(\Omega M) \simeq E^2_{0,*}(1) \ar[d]^{E^2(\tilde{i}_*)} & = E^3_{0,*}(1) & =
...& = E^{\infty}_{0,*}(1) \ar[d]^{E^{\infty}(\tilde{i})} & =
H_*(\Omega M) \ar[d]^{\tilde{i}}   \\
H_*(\Omega M) \simeq \mathbb{E}^2_{-d,*}(2) \ar@{->>}[r]&
\mathbb{E}^3_{0,*}(2) \ar@{->>}[r] &  ... \ar@{->>}[r]& 
\mathbb{E}^{\infty}_{0,*}(2) \ar@{^{(}->}[r] & \mathbb{H}_*(LM)}. 
$$
Since $\tilde{i}_*$ is injective and $H_k(\Omega M)$, $H_k(LM)$
are of finite type, the surjective maps at the bottom of the
diagram are in fact equalities. This proves that there is no non-zero
differentials arriving on $\mathbb{E}^*_{-d,*}(2)$. We conclude with
the above lemma that all the differentials of the Cohen-Jones-Yan spectral
sequence are zero.

\rightline{$\square$}

{\bf 4.6 Spheres.}
{\it Assume $n \geq 2$.

$I:\mathbb{H}_k(LS^{2n-1}) \to H_k(\Omega S^{2n-1})$ is an isomorphism
for $k = 2ni$ $i \geq 0$, $0$ otherwise.

$I:\mathbb{H}_k(LS^{2n}) \to H_k(\Omega S^{2n})$ is an isomorphism
for $k = 2i(2n-1)$ $i \geq 0$, $0$ otherwise.}

Proof:  
In \cite{CJY},  Cohen Jones and Yan have shown that all the
differentials of their spectral sequence are zero for odd spheres. Then Theorem C proves that
$im(I)=H_*(\Omega S^{2n+1})$.

For the case of even dimensionnal spheres (except the $2$-sphere), we
need the results of \cite{FHV} wich proves the following result
with rational coefficients:
$im(I)= H_k(\Omega S^{2n}; \mathbb{Q}) $ $k = 2i(n-1)$, $0$ elsewhere.
This result gives the image of the torsion free part of $I$. For
degree reasons, the image of the torsion part is zero. 

\rightline{$\square$}

{\bf 4.7 Stiefel manifolds.}
Consider the fibration $S^5 \to SO(7)/SO(5) \to S^6$.
Using together Theorems A and B, we prove that the differentials of the  Cohen-Jones-Yan
spectral sequence are zero until level $2$ while the extension issues are not
trivial. Applying theorem C, we
obtain that $im(I)=H_*(\Omega (SO(7)/SO(5))) = \mathbb{Z}[a] \otimes
\mathbb{Z}_2[b]$ with $deg(a) = 2(6-1) = 10$ and $deg(b)=5-1=4$.

\vspace{3mm}
\centerline{\sc 5. Application of the main result to the space of free paths.}

{\bf 5.1}
The last application uses the composition product on the space of free paths of
$M$, denoted by $M^I$. More explicitely, for two paths
$\gamma_1,\gamma_2 \in M^I$ such that $\gamma_1(1) = \gamma_2(0)$,
we denote by $\gamma_1 * \gamma_2$  the composed path. On homology, we define
the {\it path product} $\tilde{\mu}=Comp_* \circ \hat{i}_!: H_*(M^I) \otimes H_*(M^I) \longrightarrow
H_{*-m}(M^I)$ from the following diagram:

$\xymatrix{
M^I \times M^I  & M^I \times_{M^3} M^I \ar[l]_{\hat{i}} \ar[r]^{Comp} & M^I}$
where $\hat{i}$ is the canonical inclusion.

{\bf 5.2 Proposition 4}
{\it  
The path product on $\mathbb{H}_*(M^I)$ is identified, via
the isomorphism $H_*(M^I) \simeq H_*(M)$ to the intersection product
on $\mathbb{H}_*(M)$.}

Proof:
The homotopy $\mathcal{H}: LM \times I \longrightarrow M^I \qquad
\gamma ,s \longmapsto (t \mapsto \gamma(st) )$ proves that the
inclusion $j: LM \hookrightarrow M^I$ is homotopic to the projection
$ev(0): LM \longrightarrow M $. Since composition of paths on $M^I$
restricted to $LM$ is the composition of loops, the following diagram
is commutative:
$$
\xymatrix{
\mathbb{H}_*(M^I) \otimes  \mathbb{H}_*(M^I) \ar[d]^{\tilde{\mu}} \ar@{<->}[r] &
\mathbb{H}_*(M) \otimes  \mathbb{H}_*(M) \ar[d]_{\bullet} & &
\mathbb{H}_*(LM) \otimes \mathbb{H}_*(LM) \ar[ll]_{ev(0) \otimes ev(0)} \ar[d]^{\circ} \\
\mathbb{H}_*(M^I) \ar@{<->}[r] & \mathbb{H}_*(M) & & \mathbb{H}_*(M^I)
\ar[ll]_{ev(0)}}
$$
where $\circ$ is the Chas and Sullivan loop product, $\bullet$ is the
intersection product and $\tilde{\mu}$ is the path product on
$\mathbb{H}_*(M^I)$. This ends the proof of Proposition 2.

\rightline{$\square$}

{\bf 5.3}
Denote $\bold{H}_*(M \times M ) = H_{*+d}(M \times M) =
\mathbb{H}_{*-d}(M \times M)$. Define the new
product 
$$ \diamond : \bh_*(M \times M ) \otimes \bh_*(M \times M )
\longrightarrow \bh_*(M \times M )$$
$$(a \times b) \otimes (c \times d) \longmapsto p_{2*}(a \times (b
\bullet c) \times d)$$
where $p_2: M \times M \times M \longrightarrow M
\times M$ is the projection on the first and the third factor of $ M \times M
\times M $. This product is associative, not commutative without unit
(cf example at the end).

{\bf 5.4 Theorem A''}
{\it  Let $M$ be a smooth closed $m$-dimensional
oriented manifold. There is a multiplicative structure on the
$(d,0)$-regraded Serre spectral sequence associated to the fibration 
$
\xymatrix{
\Omega M \ar[rr] & & M^I \ar[rr]^{(ev_1,ev_0)} & &M \times M.}
$
Furthermore, if we suppose that $\pi_1(M)$ acts trivially on $\Omega M$
then we have at the $E^2$-level: 

$\be^2_{*,*} = \bh_*(M \times M; H_*(\Omega
M))$ contains $\bh_*(M \times M) \otimes  H_*(\Omega M)$ as
subalgebra. The structure of algebra on $\bh_*(M \times M) \otimes
H_*(\Omega M)$ is given by $\diamond$ on $\bh_*(M \times M)$ and by
the Pontryagin product on $H_*(\Omega M)$. Furthermore, we have $\be^2_{*,*} \Rightarrow \mathbb{H}_*(M)$ as algebra for the
intersection product.}

Proof: Consider the following commutative diagram:
$$
\xymatrix{
\Omega M \times \Omega M \ar[rr] & & M^I \times M^I \ar[rrr]^{(ev(0),
ev(1)) \times (ev(0), ev(1))} & & & M \times M \times M \times M \\
\Omega M \times \Omega M \ar[rr] \ar[u]^{id} \ar[d]^{Comp}& & M^I \times_{M^3} M^I
\ar[rrr] \ar@{^{(}->}[u]^{\tilde{D}} \ar[d]^{Comp} & & & M \times M \times M  \ar@{^{(}->}[u]^D \ar[d]^{p_2}\\
\Omega M \ar[rr] & & M^I  \ar[rrr]^{(ev(0),
ev(1))} & & & M \times M  }
$$
where $D:M^3 \longrightarrow M^4 \quad , \quad (x,y,z) \longmapsto
(x,y,y,z)$ and $\tilde{D}$ is defined from $D$ by pull-back. Let us
denote $\gamma$ the composition of paths or pointed loops.
The upper part of the diagram is a fiber embedding, so that, applying
the main result, $\tilde{D}_!$ induces the morphism of Serre spectral
sequences $E^*_{*,*}(\tilde{D}_!)$. The lower part of the diagram is a
morphism of fibration. Then $E^*_{*,*}(Comp_* \circ \tilde{D}_!):
\be^*_{*,*}(M^I \times M^I) \longrightarrow \be^*_{*,*}(M^I)$ 
is a morphism of spectral sequence witch provides the annonced
multiplicativity of the Serre spectral sequence associated to the
fibration $\Omega M \longrightarrow M^I \longrightarrow M \times M $.

\rightline{$\square$} 

{\bf 5.5 Remark} In the case of a fibred space, we could state the
  Theorem B'' analogous of Theorems B and B' , but it is not necessary since, by Proposition 4, this
  Theorem B'' is in fact Proposition 1.

{\bf 5.6 Example: the $\diamond$ product on $\bh_*(S^3 \times S^3)$.}

We apply the above result to the
fibration
$
\xymatrix{
\Omega S^3 \ar[rr] & & {S^3}^I \ar[rr]^{(ev(1),ev(0))} & &S^3 \times S^3.}
$
Denote by $[S^3] \in H_3(S^3)$ the fundamental class and by
$1$ a generator of $H_0(S^3)$. Then, we obtain the following ''table
of multiplication'' for $\diamond$ :
$$
\begin{array}{c|c|c|c|c|}
\diamond & 1 \times 1 & [S^3] \times 1 & 1 \times [S^3] & [S^3] \times
[S^3] \\
\hline\\
 1 \times 1 & 0 & 1 \times 1 & 0 &  1 \times [S^3] \\
\hline\\
\lbrack S^3 \rbrack \times 1 & 0 & [S^3] \times 1 & 0 & [S^3] \times
[S^3] \\
\hline\\
 1 \times [S^3] & 1 \times 1 & 0 &  1 \times [S^3] & 0 \\
\hline\\
\lbrack S^3 \rbrack \times \lbrack S^3 \rbrack  &  1 \times [S^3] & 0
&  [S^3] \times [S^3] & 0 \\
\hline\\
\end{array}
$$

Here we have: $\be^2_{*,*} = \be^2_{-3,*} \oplus  \be^2_{0,*} \oplus  \be^2_{3,*}$
We denote by $1_{\Omega}$ a generator of $H_0(\Omega S^3)$ and by
$u$ a generator of  $H_2(\Omega S^3)$. We put $a = ( 1 \times [S^3] +
[S^3] \times 1) \otimes 1_{\Omega}$, $b = ( 1 \times [S^3] - [S^3]
\times 1) \otimes 1_{\Omega}$ and $c = ([S^3] \times [S^3]) \otimes 1_{\Omega}$ . The only non zero differential is $d_3$
and we have $d_3(a) = 0$, $d_3(b) = (1 \times 1) \otimes 1_{\Omega}$
and $d_3(c) \neq 0 $ lies in $\be^3_{0,2}$. At the aboutment,
it remains only $(1 \times 1) \otimes 1_{\Omega} $ representing $1$ in
$\mathbb{H}_{-3}(S^3)$ and $a$ representing $[S^3]$ in
$\mathbb{H}_0(S^3)$. Let us denote by $\circ$ the induced product on
the shifted spectral sequence. We check that $ a \circ a =  ( 1 \times [S^3] +
[S^3] \times 1) \otimes 1_{\Omega} \quad \circ \quad ( 1 \times [S^3] +
[S^3] \times 1) \otimes 1_{\Omega} =  ( 1 \times [S^3] +
[S^3] \times 1) \diamond ( 1 \times [S^3] + [S^3] \times 1)  \otimes
(1_{\Omega}* 1_{\Omega}) =  ( 1 \times [S^3] +
[S^3] \times 1) \otimes 1_{\Omega} = a $ wich correspond to $[S^3]
\bullet [S^3] = [S^3]$ in $\mathbb{H}_*(S^3)$. In the same way, we check that $(1 \times 1) \otimes
1_{\Omega} \quad \circ \quad (1 \times 1) \otimes 1_{\Omega} = 0 $ and
that $(1 \times 1) \otimes 1_{\Omega} \quad \circ \quad ( 1 \times [S^3] +
[S^3] \times 1) \otimes 1_{\Omega} = ( 1 \times [S^3] +
[S^3] \times 1) \otimes 1_{\Omega} $ and $( 1 \times [S^3] + [S^3]
\times 1) \otimes 1_{\Omega}  \quad \circ \quad (1 \times 1) \otimes
1_{\Omega} =  ( 1 \times [S^3] + [S^3] \times 1) \otimes 1_{\Omega} $.
Thus we recover the intersection product on $\mathbb{H}_*(S^3)$.

{\bf 5.7 Remark}
As a final remark, let us consider the fiber embedding:
$$
\xymatrix{
\Omega M  \ar[d] \ar[r]^{id} & \Omega M \ar[d] \\
 LM   \ar[d]^{ev(0)} \ar@{^{(}->}[r]^{\tilde{\Delta}} & M^I
 \ar[d]^{ev(0) \times ev(1)}\\ 
M \ar@{^{(}->}[r]^{\Delta}  & M \times M } 
$$
($\Delta$ is the diagonal embedding).
Applying the main result, there is a morphism of spectral sequences:
$E(\tilde{\Delta}_!):\be^*_{*,*} \longrightarrow \mathbb{E}^*_{*-d,*}$
given at the $E^2$-level by:
$E^2(\tilde{\Delta}_!):\be^2_{*,*}=\bh_*(M \times M ;
\mathcal{H}_*(\Omega M)) \longrightarrow \mathbb{H}_{*-d}(M ;
\mathcal{H}_*(\Omega M)) \,, \quad (x \times y ; \omega) \longmapsto (x \bullet y ; \omega).$
This morphism of spectral sequences is not multiplicative.

\end{document}